\documentclass[11pt]{amsart}
\usepackage{amssymb}
\newtheorem{lemma}{\bf Lemma}[section]
\newtheorem{theorem}[lemma]{\bf Theorem}
\newtheorem{corollary}[lemma]{\bf Corollary}
\newtheorem{proposition}[lemma]{\bf Proposition}
\newtheorem{example}[lemma]{\bf Example}
\newtheorem{definition}[lemma]{\bf Definition}
\newtheorem{problem}[lemma]{\bf Question}
\begin{document}
\parskip = 0mm
\title{On the closure of the diagonal of a $T_1$-space}
\author{Maria Luisa Colasante}
\author{Dominic van der Zypen}

\address[M.L.Colasante]{Departamento de Matem\'{a}ticas,
Facultad de Ciencias,
Universidad de Los Andes,
M\'{e}rida 5101, Venezuela,
{\tt marucola@ula.ve}}

\address[D.van der Zypen]{Allianz Suisse Insurance Company,
Laupenstrasse 27,
CH-3001 Bern, Switzerland,
{\tt dominic.zypen@gmail.com}}



\begin{abstract}
Let $X$ be a topological space. The 
closure of $\Delta = \{(x,x):x\in X\}$ in $X\times X$ is a symmetric
relation on $X$. We characterise those equivalence relations on an infinite set 
that arise as the closure of the diagonal with respect to a $T_1$-
topology. 
\end{abstract}

\maketitle

\def\Qa{\mathbb{Q}_0}
\def\Qb{\mathbb{Q}_1}
\def\Q{\mathbb{Q}}
\def\card{{\rm card}}
\parskip = 2mm
\parindent = 0mm
\def\Part{{\rm Part}}
\def\P{{\mathcal P}}
\def\Eq{{\rm Eq}}
\def\cld{Cl_\tau(\Delta)}
\def\Csing{{\mathcal C}_{\{*\}}}
\def\Cftwo{{\mathcal C}_{{\rm fin}>1}}
\def\Cinf{{\mathcal C}_{\infty}}
\def\Pcf{{\mathcal P}_{\rm cf}}
\def\Fn{{\mathcal F}_n}
\def\proof{{\it Proof. }}


\vspace*{-4mm} 

\section{Introduction}

Our starting point is the following well-known proposition. A topological
space $(X,\tau)$ is $T_2$ if and only if the diagonal $\Delta =
\{(x,x): x\in X\}$ is a closed subset of $X\times X$.
\parskip = 2mm

Generally speaking, the closure of the diagonal
is a reflexive relation on $X$. That relation can be used to characterise
further separation axioms on topological spaces. This has been done in 
\cite{Atlanta}. 

We use the following simple observation throughout the paper:

\begin{quote} \label{obs} 
Let $(X,\tau)$ be a topological space and let $x,y\in X$. Then
$(x,y) \notin \cld$ if and only if $x$ and $y$ have disjoint open neighborhoods.
\end{quote}

This statement is straightforward to prove; moreover it implies that $\cld$ is
always a symmetric relation. We make use of the observation without
explicitly mentioning it.

As a simple example, consider any infinite set $X$ and equip it with 
the topology $\Pcf(X)$, ie.~the collection of all cofinite sets together
with the empty set. 
Note that a topology $\tau$ is $T_1$ iff $\tau \supseteq
\Pcf(X)$. In $\Pcf(X)$ any two members intersect. So by the observation
made above one obtains $\cld = X\times X$.

Moreover note that for two topologies $\sigma \supseteq \tau$ one obtains 
${Cl_\sigma(\Delta)} \subseteq \cld$. Which reflexive and symmetric
relations on a set $X$ can be represented as the closure of the 
diagonal on $X\times X$ with respect to some topology? This question
is very natural and is still open. We confine ourselves
to equivalence relations and $T_1$-spaces. 

\begin{definition}
Let $R$ be an equivalence relation on an infinite set $X$. Then 
$R$ is said to be {\bf $T_1$-realisable} if there is a $T_1$ topology
on $X$ such that $R = \cld$.
\end{definition}
The main goal of this article is to characterise $T_1$-realisable 
equivalence relations and we provide a characterisation in theorem \ref{charac}.
Moreover, we show in Example \ref{theExample} that the closure
of the diagonal need not be transitive.

The following tool is straightforward to prove.
It allows us to consider just
{\it bases} instead of whole topologies when dealing with the closure
of the diagonal. By a basis we mean a collection $\mathcal B$ of open subsets
of a set $X$, such that for all $B_1,B_2 \in {\mathcal B}$ one obtains
\begin{quote}
for every $x\in B_1\cap B_2$ there is $B\in {\mathcal B}$ such that 
$x\in B$ and $B\subseteq B_1\cap B_2$.
\end{quote}
\begin{lemma} \label{baslem}
Let $X$ be a topological space and let ${\mathcal B}$ be a basis
of $X$. Then $(x,y)\notin \cld$ if and only if there exist disjoint members
$A, B\in {\mathcal B}$ such that $x\in A$ and $y\in B$.
\end{lemma}
If we consider all topologies on a set $X$, then every
equivalence relation can be written as the closure of the diagonal
of some topology, as the following shows:

Let $R$ be an equivalence relation on $X$
For each $x\in X,$ we let $R(x)=\{y\in X: (x,y)\in R\}$ denote the equivalence
class of $x$. Let
\begin{center}
$\tau _{R}=\{V\subset X:$ $x\in V$ implies $R(x)$ $\subset V\}.$
\end{center}
It is easy to see that $\{R(x): x\in X\}$ is a basis for the topology
and that therefore $x,y$ can be separated by disjoint open sets iff 
$(x,y)\notin R$. 

Next we introduce some notation. It is well known that equivalence relations
and partitions of a set $X$ are in a natural correspondence. For
an equivalence relation $R$ let $\Part (R)=\{R(x):x\in X\}$ and each
partition $\P$ on $X$ let $\Eq (\P) = \{(x,y)\in X\times X: (\exists B \in \P):
x,y\in P\}$. It is obvious that $\Part (R)$ is a partition and $\Eq(\P)$
is an equivalence relation on $X$. The elements of a partition are 
called ``blocks''. We also refer to the the equivalence classes of an
equivalence relation $R$ (ie the blocks of $\Part (R)$) as ``blocks''.

\section{Infinite blocks only}
We present a solution for the following particular case: $R$ is an equivalence
relation on a set $X$ such that $R(x)$ is infinite
for each $x\in X.$
\begin{proposition}\label{infblocks}
If $R$ is an equivalence relation on an infinite set $X$ such that every
block of $R$ is infinite, then $R$ is $T_1$-realisable.
\end{proposition}
\proof
Consider the following topology on $X$:
\begin{center}
$\tau =\{X,\emptyset \}\cup \{V\subset X:$ $x\in V$ implies $R(x) \cup V
$ is cofinite in$ R(x)\}$
\end{center}
It is easy to check that $\tau $ is $T_{1}$.
Suppose that $(x,y)\in R$ and pick any open neighborhoods $U,V$ of $x$ 
and $y$, respectively. Clearly, $U\cap R(x)$ and $V\cap R(x)$ are cofinite
subsets of the infinite set $R(x)$ by definition of $\tau$, so they intersect.
Conversely suppose $(x,y)\notin R$. Then $R(x)$ and $R(y)$ are disjoint
open neighborhoods of $x$ and $y$, respectively. So
one concludes that $R=\cld$. \qed

\section{Allowing for finite blocks}
For any equivalence relation $R$ on a set $X$ we define three important sets:
$\Csing(R)$, the set of all singleton blocks, $\Cftwo(R)$, the set of all finite
blocks containing more than one element, and $\Cinf(R)$, the 
set of all infinite blocks.

More formally, 
we let $\Csing(R)= \{B\in \Part(R): B = \{x\} \textrm{ for some } x\in X\}$ 
and $\Cinf(R) =  \{B\in \Part (R): B\textrm{ is infinite}\}$ and last
$\Cftwo(R) = \Part(R) \setminus (\Csing(R) \cup \Cinf(R))$.
\begin{proposition}\label{easyprop}
If all blocks of an equivalence relation $R$ are infinite or singletons,
$R$ is $T_1$-realisable.
\end{proposition}
\proof We just give a sketch. Let $Y = \bigcup \Cinf(R)$ and let
$Z = \bigcup \Csing(R)$. Clearly, the restriction of $R$ to $Y$
has only infinite blocks. Endow $Y$ with the topology described in proposition
\ref{infblocks} and give $Z$ the discrete topology. It is easy to see
that this is a basis for a topology on $X=Y\cup Z$ such that $R=\cld$.
\qed

If we allow for finite blocks with more than one element, considerations get 
more involved. 

\begin{proposition} \label{fintwo}
Let $R$ be an equivalence relation on an infinite
set $X$ such that $\Cftwo(R)$ is finite and nonempty.
\begin{enumerate}
\item If $\Part(R)$ is finite then $R$ is not $T_1$-realisable.
\item If $\Part(R)$ is infinite then $R$ is $T_1$-realisable.
\end{enumerate}
\end{proposition}
\proof
(1) Assume that there is a $T_1$ topology $\tau$ such 
that $\cld = R$. 
Pick $a\in X$ such that $R(a)$ is finite and there exists
$b\neq a$ such that $(a,b)\in R$. 

{\it Claim 1.} Every open neighborhood of $a$ is infinite. - If there
were a finite neighborhood of $a$ then $\{a\}$ would be open since $\tau$
is $T_1$. So $\{a\}$ and $X\setminus \{a\}$ separate $a$ and $b$, 
which implies $(a,b) \notin \cld = R$, contradicting our choice of
$a,b$.

{\it Claim 2.} There exists $B^* \in \Cinf(R)$ such that
for every open neighborhood $U$ of $a$ one obtains $B^*\cap U \neq \emptyset$. 
Suppose the contrary, so for all $B\in \Cinf(R)$ there is an open
neighborhood $U_B$ of $a$ such that $U_B\cap B = \emptyset$. Then
$U' = \bigcap\{U_B: B\in \Cinf(R)\}$ is open since $\Cinf(R)$ is finite,
and $U' \subseteq \bigcup \Cftwo(R) \cup \bigcup \Csing(R)$, which
is a finite set. So $U'$ is a finite neighborhood of $a$, contradicting
claim 1. So claim 2 is proved.

Now pick $B^*$ from claim 2 and let $x\in B^*$. We show that
$x$ and $a$ cannot be separated by open sets, which then contradicts
$(a,x) \notin R$. - Pick open neighborhoods $U,V$ of $a$ and $x$, respectively.
Using claim 2, pick $y\in B^* \cap U$. Since $(x,y) \in R$, they cannot
be separated by open neighborhoods and since $U$ is an open neighborhood
of $y$, one obtains that $U$ and $V$ intersect, and we are done.

(2) We distinguish two cases:

{\sl Case 1.} $\Csing(R)$ is infinite. Recall that $\Cftwo(R)$ is finite
by assumption of the proposition. For each $C\in \Cftwo(R)$ pick
an infinite subset $S_C\subseteq \bigcup \Csing(R)$ such that if 
$C,D\in \Cftwo(R)$
then $S_C\cap S_D = \emptyset$. 
We give a basis for a topology $\tau$ by
\begin{tabbing}
pr\= bs \= eq \= some more some more \= endend \kill
\> ${\mathcal B}$ \> $=$ \> $\{\{x\}: x\in \bigcup\Csing(R)\} \hspace*{2mm} \cup$ \> \\
\> \> \> $\{U\subseteq X: U\in \Pcf(B) \textrm{ for some }B\in \Cinf(R)\} 
\hspace*{2mm} \cup$ \> \\
\> \> \> $\{V\subseteq X: V\cap (\bigcup\Cftwo(R)) = \{x\} \textrm{ for some
} x \textrm{ and } V = \{x\} \cup U$ \> \\
\> \> \> \> $\textrm{ for some } U\in \Pcf(S_{R(x)})\}$.\\
\end{tabbing}
We argue shortly that ${\mathcal B}$ is indeed a basis. Designate the three
``parts'' of ${\mathcal B}$ by $P_1, P_2, P_3$ such that ${\mathcal B}$ is the
disjoint union of $P_1, P_2$ and $P_3$. Clearly, the intersection of
two members of $P_i$ for $i=1,2$ is empty or again in $P_i$. Let $V,W\in P_3$
and suppose that $V\cap (\bigcup\Cftwo(R)) = \{x\}$ and 
$W\cap (\bigcup\Cftwo(R)) = \{y\}$. Moreover, let $V = \{x\} \cup
U$ for some $U\in \Pcf(S_{R(x)})$ and let $W = \{y\} \cup
U'$ for some $U'\in \Pcf(S_{R(y)})$. If $x\neq y$ then $V\cap W = \emptyset$.
If $x=y$ then $U\cap U'\in \Pcf(S_{R(x)})$ and $V\cap W = \{x\} \cup (
U\cap U')$. So $V\cap W \in {\mathcal B}$. Last, let $M_i\in P_i$. We get $M_1\cap M_2 = M_2 \cap M_3
= \emptyset$ and $M_1\cap M_3$ is empty or contains for each
$x\in M_1\cap M_3$ a member $B$ of ${\mathcal B}$ such that $x\in B\subseteq M_1\cap
M_3$: take $B=\{x\}$.

Moreover, a case
distinction shows that $x,y\in X$ can be separated by disjoint members
of $\mathcal B$ if and only if $(x,y)\notin R$. And for $x\neq y$ in $X$ there
are basic open sets containing $x$ but not $y$ and vice versa. So the
topology generated by ${\mathcal B}$ is $T_1$.

{\sl Case 2.} $\Csing(R)$ is finite, so $\Cinf(R)$ is infinite
since $\Part(R)$ is infinite. Recall that $\Cftwo(R)$ is finite
by assumption of the proposition. For each $C\in \Cftwo(R)$ pick
an infinite subset ${\mathcal S}_C\subseteq \Cinf(R)$ such that if $C,D\in \Cftwo(R)$
then ${\mathcal S}_C\cap {\mathcal S}_D = \emptyset$. Note that there is 
a subtle difference in
the definition of $S_C$ in the above case and ${\mathcal S}_C$ here: 
in case 1, $S_C$ was 
a subset of $X$, and here ${\mathcal S}_C$ is a subset of $\Cinf(R)$.
We give a basis for a topology $\tau$ by
\begin{tabbing}
pr\= bs \= eq \= some more some more \= endend \kill
\> ${\mathcal B}$ \> $=$ \> $\{\{x\}: x\in \bigcup\Csing(R)\} \hspace*{2mm} \cup$ \> \\
\> \> \> $\{U\subseteq X: U\in \Pcf(B) \textrm{ for some }B\in \Cinf(R)\} 
\hspace*{2mm} \cup $\> \\
\> \> \> $\{V\subseteq X: V\cap (\bigcup\Cftwo(R)) = \{x\} \textrm{ for some
} x \textrm{ and } V = \{x\} \cup (\bigcup{\mathcal T})$ \> \\
\> \> \> \> $\textrm{ for some } {\mathcal T} \in \Pcf({\mathcal S}_{R(x)})\}.$
\end{tabbing}
In a very similar way to what we did in case 1, we can check that
${\mathcal B}$ is a basis. And again, the topology generated by ${\mathcal B}$ is
$T_1$. Moreover, a case
distinction shows that $x,y\in X$ can be separated by disjoint members
of $\mathcal B$ if and only if $(x,y)\notin R$.
\qed

\begin{corollary}
On an infinite set, an equivalence relation with finitely many blocks
is $T_1$-realisable if and only if every finite block is a singleton.
\end{corollary}
Proposition \ref{fintwo} describes what happens if $\Cftwo(R)$ is finite
and nonempty. Next we we look at what happens if $\Cftwo(R)$ is infinite. 
First we prove that if each block consists of exactly two points, then
the relation is $T_1$-realisable.
\begin{lemma}\label{lem2} 
Let $X$ be an infinite set and $R$ be an equivalence relation
such that $\card(R(x))=2$ for all $x\in X$. Then $R$ is $T_1$-realisable.
\end{lemma}
\proof
We construct a topological space $(Z,\sigma)$ whose ground set
$Z$ is equinumerous to $X$ and equip
it with an equivalence relation $S$ such 
that $\card(S(z)) = 2$ for all $z\in Z$ such
that the above-mentioned properties are satisfied.
It is not difficult to prove that there is a bijection $\varphi: Z\to X$
such that $(z_1, z_2)\in S$ if and only if $(\varphi(z_1),\varphi(z_2))\in R$.
Then $\tau = \{\varphi(U): U\in \sigma\}$ is a topology satisfying the
condition of the lemma.
Let $Z=X\times\Q\times\{0,1\}$. Note that clearly $Z$ is 
equinumerous with $X$. Moreover we let the equivalence relation $S$
on $Z$ be defined by
$$(x,q,k) \sim_S (x',q',k') \textrm{ iff } x=x' \textrm{ and } q=q'.$$
Clearly, each block of $S$ has $2$ elements.

We define $B\subseteq 
Z$ as basic open if and only if there are $x\in X, q\in \Q$ and $\delta \in \Q_{>0}$
such that $B$ is a cofinite subset of 
$$B_\delta(x,q):=\{x\}\hspace*{1mm} \times\hspace*{1mm}
B_\delta(q)\hspace*{1mm} \times \hspace*{1mm}
\{0,1\},$$
where $B_\delta(q)=\{q'\in \Q: |q'-q|<\delta\}$. Intuitively speaking,
$B_\delta(x,q)$ is a series of $2$ copies of $B_\delta(q)$, such that 
each copy lies at its appropriate place in the set
$X\times\Q\times \{k\}$ for each $k$.

It is easy to verify that the collection of all basic open sets
is a basis that gives rise to a $T_1$-topology. Moreover the following
is readily verified: members $(x,q,k)$ and $(x',q',k')$ of $Z$
can be separated by basic open sets if and only if $x\neq x'$ or $q\neq q'$.
So using lemma \ref{baslem} and the reasoning at the beginning
of this proof we are done.
\qed
\begin{proposition} \label{Cftwoprop}
Let $X$ be an infinite set and $R$ an equivalence
relation on $X$. Suppose that $\Part(R)$ is infinite and
each block has more than one element.
Then $R$ is $T_1$-realisable.
\end{proposition}
\proof
For each $B\in \Part(R)$ pick two distinct representatives 
$r_1(B), r_2(B)\in B$.
Let $W=\{r_1(B): B\in \Part(R)\} \cup \{r_2(B): B\in \Part(R)\}$ and let
$S=R\cap (W\times W)$. Clearly, each block of $S$ has just two 
elements, so by lemma \ref{lem2} there is a $T_1$-topology $\sigma$ on $W$
such that $Cl_\sigma(\Delta_W)=S$. 

Using $\sigma$ we equip $X$ with a topology having the desired
property. We say that $V\subseteq X$ is basic open if and only 
if one of the following two conditions holds:
\begin{enumerate}
\item $V\in \sigma$, that is $V\subseteq W$ and $V$ is open;
\item there is $x\in X\setminus W$ and $A\in \sigma$
with $A\cup \{r_1(R(x))\}\in \sigma$ and $B=\{x\}\cup A$.
\end{enumerate}
As we easily verify, the collection ${\mathcal B}$ of
basic open elements is indeed a basis and very intuitively speaking
``things in the basis happen more or less on $W$''. We define
$\tau$ to be the topology generated by $\mathcal B$.

Suppose that $(x,y)\in R$ and that $V_x, V_y$ are basic open sets
containing $x$ and $y$, respectively. So this leads to case distinction.
Suppose that $x,y\notin \{r_1(R(x)), r_2(R(x))\}$.
Then there are $A_x, A_y$ in $\sigma$
such that $V_x =\{x\} \cup A_x$ and $V_y=\{y\}\cup A_y$ and 
$A_x\cup\{r_1(R(x))\}\in \sigma$ and
$A_y\cup\{r_1(R(y))\}\in \sigma$ where
of course $R(x)=R(y)$. So $(A_x\cap A_y)\cup
\{r_1(R(x))\}\in \sigma$. Moreover, every neighborhood
of $r_1(R(x))$ in $(W,\sigma)$
is infinite. (Suppose otherwise: since $\sigma$ is $T_1$, we could
separate $r_1(R(x))$ and $r_2(R(x))$ by disjoint open sets 
in $W$, contradicting $(r_1(R(x)),r_2(R(x)))\in S=Cl_{\sigma}(\Delta_W)$.)
So since $(A_x\cap A_y)\cup
\{r_1(R(x))\}$ is an open neighborhood of $r_1(R(x))$, it
can not just consist of $r_1(R(x))$,
therefore $A_x, A_y$ have nonempty intersection,
and so have their respective supersets $V_x$ and $V_y$. The other cases
are treated in a similar way.

Suppose that conversely one obtains $(x,y)\notin R$. One has 
distinguish some cases. Assume that $x\notin\{r_1(R(x)), r_2(R(x))\}$ and
$y\notin\{r_1(R(y)), r_2(R(y))\}$.
Look at $r_1(R(x))$
and $r_1(R(y))$. Since $(r_1(R(x)),r_1(R(y)))
\notin S=Cl_{\sigma}(\Delta_W)$, they can be separated by
disjoint open neighborhoods $U_1, U_2$, respectively. 
By definition
of ${\mathcal B}$, the sets $\{x\}\cup U_1$ and $\{y\}\cup U_2$ are
disjoint basic open sets separating $x$ and $y$. The remaining cases are
treated in a similar way, implying that $x$ and $y$ can always be separated
by basic open sets if $(x,y)\notin R$. 

So we get $R=\cld$.\qed

\begin{corollary}If $R$ is an equivalence relation on a set $X$ such that
$\Part(R)$ is infinite, then $R$ is $T_1$-realisable.
\end{corollary}
\proof
The case that $\Cftwo(R)$ is finite has been dealt with. So suppose
that $\Cftwo(R)$ is infinite. Note that
$X$ is the disjoint union of $A=\bigcup \Csing(R)$, $B=\bigcup \Cftwo(R)$
and $C=\bigcup \Cinf(R)$. Endow $B$ with the topology described
in proposition \ref{Cftwoprop} and give $A\cup C$ the topology
constructed in \ref{easyprop}. The disjoint union of the two
topological spaces just mentioned give a $T_1$-topology on $X$ such that
$\cld = R$.
\qed

The above results can be summarised in the following theorem.
\begin{theorem}\label{charac}
Suppose that $X$ is an infinite set and $R$ is an equivalence relation
on $X$. Then the following are equivalent:
\begin{enumerate}
\item $R$ is not $T_1$-realisable;
\item $\Part(R)$ is finite and $R$ has a finite block that is not a singleton.
\end{enumerate}
\end{theorem}

\section{Further directions}
First we construct a space such that the closure of the diagonal is
not transitive.
\begin{example}\label{theExample}On $X=\omega$ we let 
\begin{itemize}
\item $D_1 = \{3n+1: n\in \omega \}$, and
\item $D_2 = \{3n+2: n\in \omega \}$.
\end{itemize}
Let $\tau$ be the topology on $\omega$ generated by the subbasis
$${\mathcal S} = \Pcf(\omega)\cup \{D_1\} \cup \{D_2\}.$$
Obviously, $\tau$ is $T_1$.

We obtain $(2,3)\in \cld$ and $(3,4)\in \cld$ but $(2,4)\notin \cld$.
So $\cld$ is not transitive.
\end{example}
In section 1
we saw that any equivalence relation is realisable by some topology,
although the topology given there is in general not even $T_0$.

\begin{proposition} \label{referee}Every equivalence relation $R$ on a set
$X$ is $T_0$-realisable.
\end{proposition}
\proof
Let $R$ be an equivalence relation on the set $X$.
For each block $B\in \Part(R)$ we pick a representative 
$r(B) \in B$ and we define $U \subseteq
X$ to be open if and only if
\begin{quote}
if $U \cap B \neq \emptyset$ for some $B\in \Part(R)$ then $r(B) \in U$.
\end{quote}
It is not difficult to see that the collection of open sets form a topology.
To see that this collection is $T_0$ pick $x\neq y \in X$. If $(x,y)\notin R$
then $R(x)$ is an open set that contains $x$ but not $y$. If $(x,y)\in R$ 
then at least one of $x,y$ does not equal $r(B)$ where $B=R(x)=R(y)$ and
we may assume $y\neq r(B)$. Then $\{x, r(B)\}$ is an open set containing
$x$ but not $y$.

Next we show that $\cld = R$. Let $(x,y) \in (X\times X) \setminus R$. 
Then $R(x)$ and $R(y)$ are disjoint open sets containing $x$ and $y$,
respectively. So $(x,y) \notin \cld$.
Conversely take $(x,y) \in R$. By construction of the open sets, every
open neighborhood of $x$ and every open neighborhood of $y$ contains
$r(R(x)) = r(R(y))$. So, $x, y$ cannot be separated by disjoint open 
sets and implies that $(x,y) \in \cld$.
\qed

The present article characterises $T_1$-realisable equivalence relations.
A natural extension of this is the examination of the realisability of
symmetric and reflexive binary relations without the requirement that
the relation be transitive. We want to conclude the article
with two open questions:
\begin{problem}
Let $X$ be any set. For which symmetric and reflexive relations $R$
does there exist a topology on $X$ such that $\cld = R$? What happens
if we confine ourselves to $T_0$- or $T_1$-spaces?
\end{problem}
The more ``natural'' spaces tend to have a transitive closure of their diagonal.
This issue is addressed in the following question:
\begin{problem} Characterise those topological spaces such that the closure
of the diagonal is transitive. Is there a ``geometrical interpretation''
of transitivity?
\end{problem}
\section{Acknowledgements}
We are grateful to both (anonymous) referees for their helpful comments.
One referee suggested the statement and the proof of proposition
\ref{referee} and pointed out an error in example \ref{theExample}. 
The other referee helped us
improve the overall structure and the language of the article.


\end{document}